\documentclass[eng]{ajceam-class2}

\usetikzlibrary{shapes.geometric, arrows}

\tikzstyle{stop} = [ellipse, minimum width=3cm, minimum height=1cm,text centered, text width=3cm, draw=black, fill=cyan!30]

\tikzstyle{start} = [ellipse, minimum width=2cm, minimum height=1cm,text centered, draw=black, fill=cyan!30]
\tikzstyle{io} = [trapezium, trapezium stretches=true, trapezium left angle=70, trapezium right angle=110, minimum width=3cm, minimum height=1cm, text centered, draw=black, fill=blue!30]

\tikzstyle{process} = [rectangle, minimum width=3cm, minimum height=1cm, text centered, text width=3cm, draw=black, fill=orange!30]

\tikzstyle{mprocess} = [rectangle, minimum width=4cm, minimum height=1cm, text centered, text width=4cm, draw=black, fill=orange!30]

\tikzstyle{lprocess} = [rectangle, minimum width=6cm, minimum height=1cm, text centered, text width=6cm, draw=black, fill=orange!30]

\tikzstyle{decision} = [diamond, minimum width=3cm, minimum height=1cm, text centered, draw=black, fill=green!30, text width=3cm]
\tikzstyle{arrow} = [thick,->,>=stealth]

\title{Reliability of Numerical Solutions in Transient Chaos}

\author[1]{Ali Goodarzi}
\author[1]{Maryam Rahimi}
\author[2]{MohammadJavad Valizadeh}
\author[3,4]{Fakhteh Ghanbarnejad}

\affil[1]{Institute of Physics, École Polytechnique Fédérale de Lausanne, Lausanne, Switzerland}
\affil[2]{Department of Mathematics, Simon Fraser University, Burnaby, Canada}
\affil[3]{Department of Physics, Sharif University of Technology, Tehran, Iran}
\affil[4]{Chair of Network Dynamics, Institute for Theoretical Physics and Center for Advancing Electronics Dresden (cfaed), Technical University of Dresden, 01062 Dresden, Germany}

\firstauthor{Goodarzi, Rahimi, Valizadeh and Ghanbarnejad}

\contactauthor{Fakhteh Ghanbarnejad}           
\email{fakhteh.ghanbarnejad@gmail.com}            

\abstract{
In dealing with nonlinear systems, it is common to use numerical solutions. Unlike the careful behavior towards the numerical results in chaotic regions, the validity of numerical results in regions of transient chaos might not always be taken into consideration. This article demonstrates that using numerical methods to solve systems undergoing transient chaos can be challenging and sometimes unreliable. To illustrate this issue, we use the Lorenz system \cite{Lorenz1} in the region of transient chaos as an example. We show how the result of the computation might completely change when using different mathematically equivalent expressions. This raises the question of which result should be relied on. To answer this question, we propose a method based on the Lyapunov exponent to determine the reliability of the numerical solution and apply it to the provided example. In fact, this method checks a necessary condition for the validity of the numerical solution. Then, by increasing the precision to the extent suggested by our method, we show that the result of our studied case passes this test. In the end, we briefly discuss the scope and limits of our method.
}

\keywords{
transient chaos, numerical solution, numerical error, Lyapunov exponent
}

\begin{document}

\maketitle
\thispagestyle{fancy}
\printcontactdata

\section{Introduction}
\label{sec:1}
\firstword{W}{hen}
it comes to studying a dynamical system, numerical methods are common, especially in many non-linear systems where an analytical solution is not an option. However, errors are inevitable in numerical approaches. Since chaos is known to be highly sensitive to any perturbation, many studies have been done to investigate the problems that these errors might cause in chaotic systems \cite{lorenz2,PhysRevLett,compnet,yao,numsol}. But what about transient chaos?
\par
Transient chaos is a common phenomenon in nonlinear dynamical systems. It can be observed in a vast variety of topics, from hydrodynamics \cite{hydrodynamics1,hydrodynamics2}, electronic circuits \cite{ecircuits1,ecircuits2}, power grid \cite{powergrid}, and NMR lasers \cite{nmr} to population dynamics \cite{popudynamics}, ecology \cite{ecology1,ecology2}, economics \cite{economics}, neural networks \cite{neuralnetworks1,neuralnetworks2} and some medical applications \cite{medapps}. Hence, it is crucial to make sure that the result of simulations in transient chaos is reliable.
\par
Transient chaos is "the form of chaos due to nonattracting chaotic sets in the phase space." \cite{lai2011transient} As the destiny of the systems undergoing transient chaos is non-chaotic, one might overlook errors in numerical solutions; because, in dynamical systems, it is usual to just focus on destiny. The detrimental effects of these errors in transient chaos are not less important than in chaos; because, in studying chaotic systems, it is common to employ a statistical perspective and describe the system based on averaging, which may diminish the effects of numerical errors. However, in transient chaos, this statistical view is not always applicable; since the errors might change the system's destiny. Also, when looking at the evolution of the system, especially in short-time transient chaos, one might mistake it for a system spiraling into a stable fixed point.
\par
This article aims to point out the effects of numerical errors on the destiny of systems undergoing transient chaos. In section \ref{sec:2}, we demonstrate this problem in a Lorenz system \cite{Lorenz1}, and after that, we discuss two possible numerical errors and the reasons that cause this problem. In section \ref{sec:4.1}, using the Lyapunov exponent, we employ a method to understand whether a numerical solution of a system is reliable. If not, we state a possible way to enhance the simulation and make it more reliable. By utilizing the introduced procedure, in section \ref{sec:4.2}, we get back to our Lorenz example and find a way to reach a reliable solution for our system.

\begin{figure*}[ht]
\begin{subfigure}{\linewidth}
    \includegraphics[width=\linewidth]{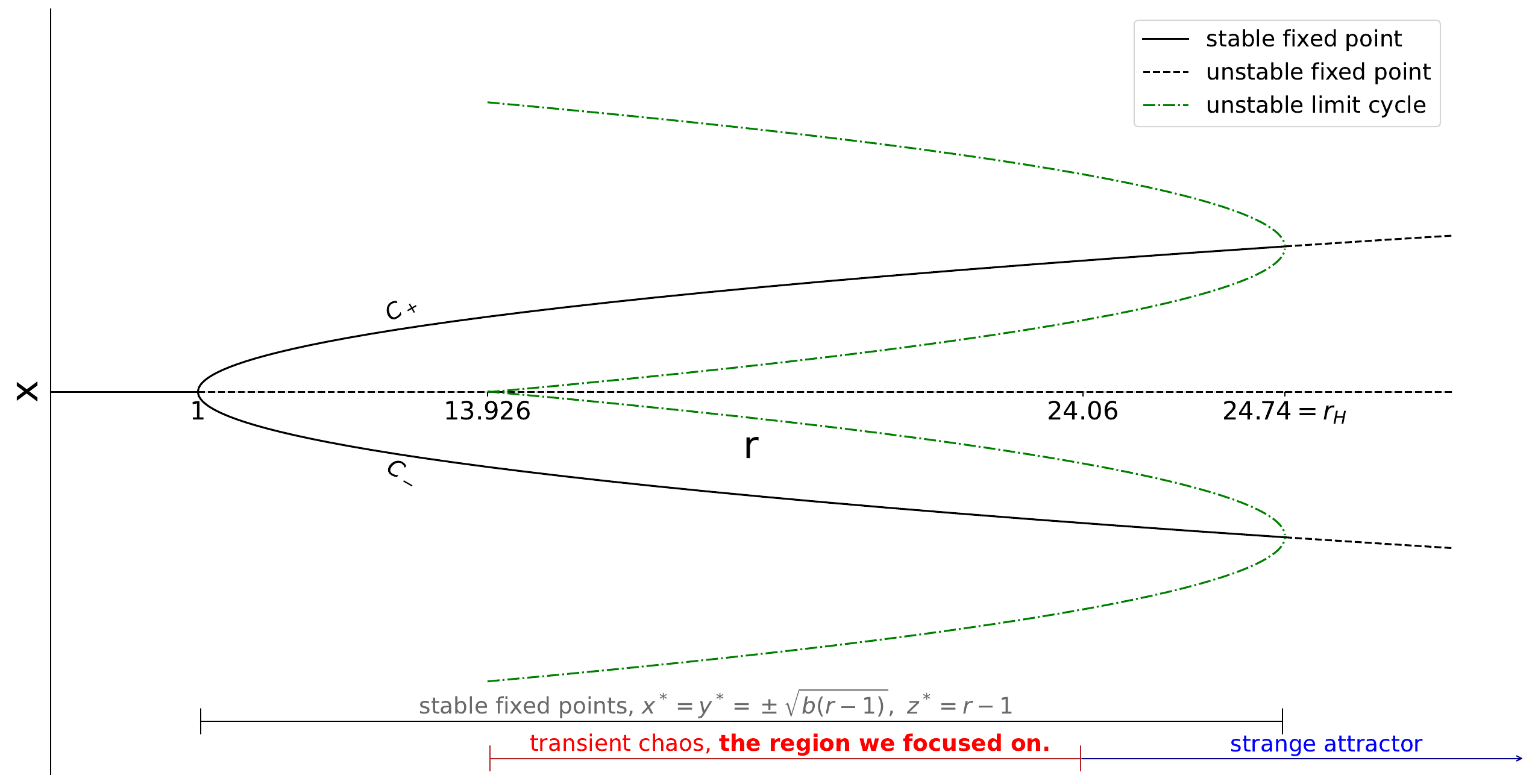}
    \end{subfigure}\notag
    \caption{
Lorenz system's behavior for small values of $r$ \cite{Strogatz2}. Note that the Lorenz equations are symmetric under $(x,y) \to (-x,-y)$. \\ In the provided example, we have focused on the region of transient chaos.}\label{fig:1}
\end{figure*}

\section{Problem Statement: which numerical result is the true destiny of the system?}
\label{sec:2}
In this article, the numerical solution is assumed valid if its results are identical with that of the theoretical solution, representing the actual outcome of the system under certain initial condition. What follows is an example to explain this problem further.
\par 
For our purpose, we use Lorenz equations \cite{Lorenz1}: 
\begin{subequations}
    \begin{align}
    \notag 
    &\dot{x} = \sigma(y-x)
    \\\label{eq:1.a} \tag{1a}
    &\dot{y} = rx-xz-y
    \\ \notag
    &\dot{z} = xy-bz
    \end{align}
\end{subequations} \noindent
According to Fig \ref{fig:1}, for $\sigma=10, b=\frac{8}{3}$, in $1<r<r_H=24.74$, the system has two stable fixed points, $C_{\pm}=(x^{\star}_{\pm},y^{\star}_{\pm},z^\star$):
\begin{equation}
    \label{eq:fixedpoints}
    x^\star_{\pm}=y^\star_{\pm}=\pm \sqrt{b(r-1)}, z^\star=r-1
\end{equation} \noindent

We concentrate on $13.926<r<24.06$, where the system can exhibit transient chaos. In this region, trajectories can wander chaotically for a while and eventually reach either of the stable fixed points.

As $r$ increases, the time that they behave chaotically grows to infinity so that in $r=24.04$, the trajectories cannot escape the strange attractor \cite{Strogatz1}.

\par

\begin{figure*}[!ht]
    \begin{subfigure}{0.5\linewidth}
    \includegraphics[width=\linewidth]{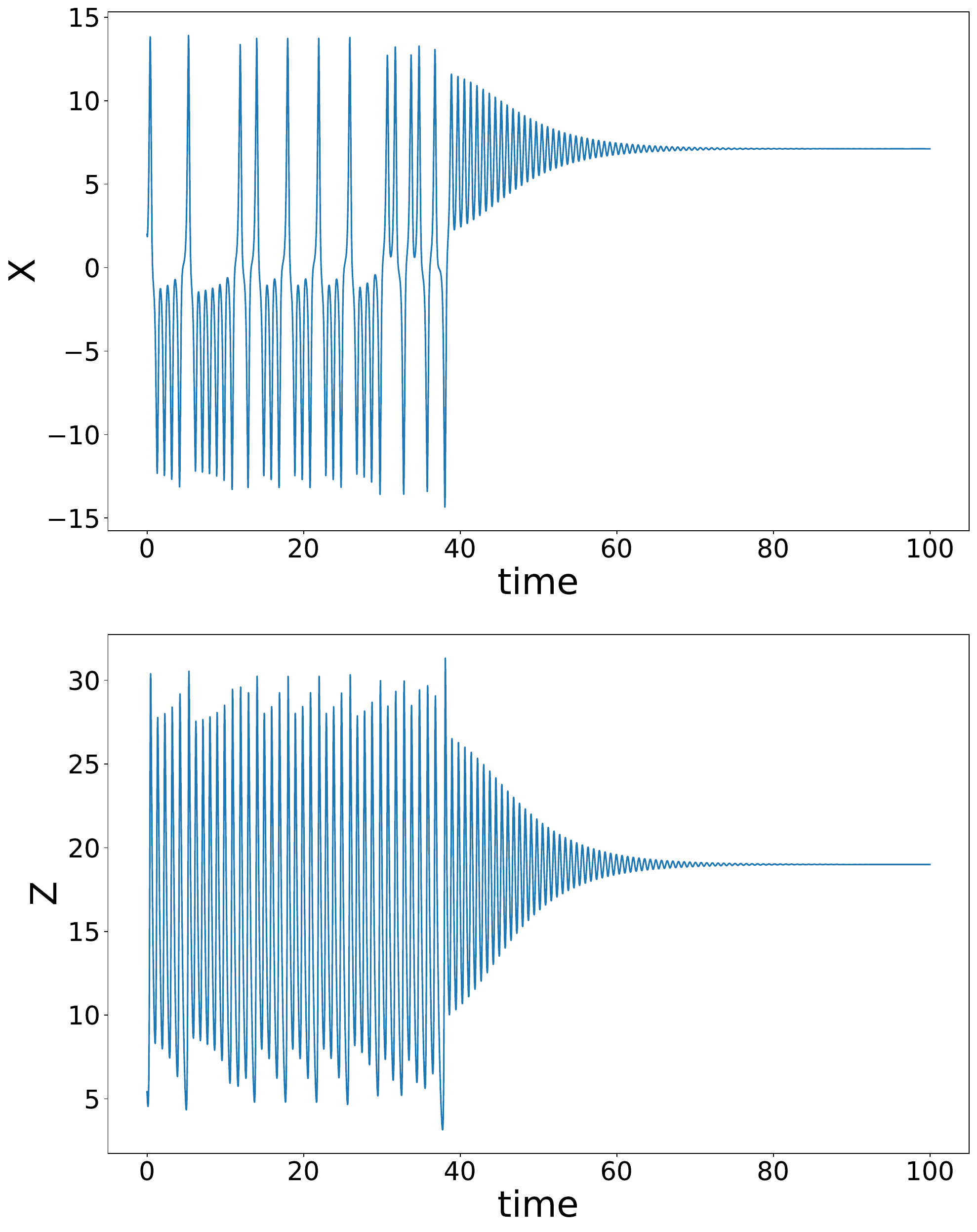} 
    \caption{with expression \eqref{eq:1.a}.}
    \label{fig:2.a}
    \end{subfigure}
    \hspace{0.2 mm}
    \begin{subfigure}{0.5\linewidth}
    \includegraphics[width=\linewidth]{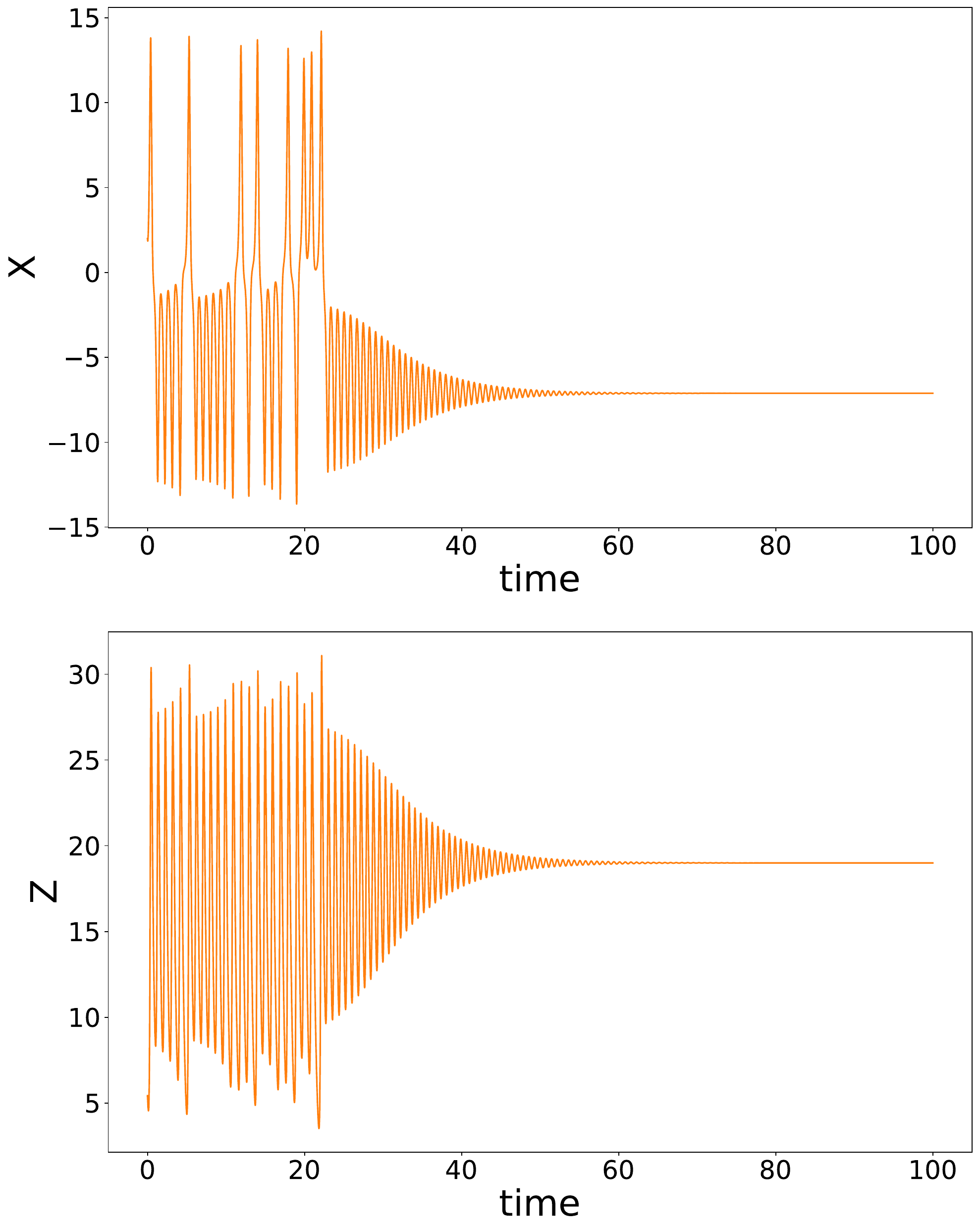}
    \caption{with expression \eqref{eq:1.b}.}
    \label{fig:2.b}
    \end{subfigure}

    \caption{Mathematically equivalent expressions reaching different fixed points. \\ Numerical results of the provided example with \textbf{NumPy float32} variables, X(t) \& Z(t) plots, where colors blue \& orange, \ref{fig:2.a} \& \ref{fig:2.b}, refer to the equations \eqref{eq:1.a} \& \eqref{eq:1.b} respectively: Lorenz equations with $r=20$ \& initial condition $(x_0,y_0,z_0)=(2,1,5.42857)$, solved with RK4 algorithm. \\Although the expressions are mathematically equivalent, the system reaches different destinations, i.e. \ref{fig:2.a} \& \ref{fig:2.b} reached fixed points $C_+$ \& $C_-$, given from equation \eqref{eq:fixedpoints},  respectively.}
    \label{fig:2}
\end{figure*}

The result of a simulation of the system for r=20, using the equation \eqref{eq:1.a}, the \emph{RK4} \footnote{\emph{Runge-Kutta 4$^{th}$ order}} algorithm \cite{kutta} \& \cite{runge} and NumPy float32 variables,  for a specific initial condition, $(x_0,y_0,z_0)=(2,1,5.42857)$, is given in Fig.\ref{fig:2.a}.
\\One can use expression \eqref{eq:1.b} instead of \eqref{eq:1.a} in the Lorenz equations for the simulation: 
\begin{equation}
    \dot{y} = rx-y-xz    \label{eq:1.b} \tag{1b}
\end{equation}

Obviously, \eqref{eq:1.a} and \eqref{eq:1.b} are mathematically equivalent, so it is convenient to expect that the simulation result should be the same. The outcome is shown in Fig.\ref{fig:2.b}.
\par 
The result is different from our expectation. Not only the numerical solution is completely different for \eqref{eq:1.a} and \eqref{eq:1.b}, but also they have reached different fixed points. 
\par
The problem is which one is the true destiny of this system? Should we rely on one of these numerical solutions or look for another way to understand the system?
\section{Types of numerical errors}
\label{sec:3}
To find out the reason for the problem, we should first be familiar with numerical errors. There are different factors in numerical approaches that cause errors in the final result:
\begin{enumerate}[label=\bf{\emph{\Alph*}})]
\item The error caused by the numerical algorithms, such as \emph{Runge-Kutta, Euler, Verlet,} etc., used to solve differential equations. Each algorithm causes \emph{local truncation} and \emph{total accumulated} errors, which depend on the size of discretization and the algorithm’s order of accuracy. For example, the \emph{RK4} algorithm has a \emph{local truncation} error of $\mathcal{O}({dt}^{5})$ and a \emph{total accumulated} error of $\mathcal{O}({dt}^{4})$, $dt$ being the discretization parameter. 

\item The other factor is that the computer uses floating-point arithmetic and has a finite precision, causing rounding errors. The order of this rounding error depends on the type of variables used in the simulation (the number of digits that can be stored) and the order of magnitude of numbers used in the operations. (The error of each operation is equal to the precision of the larger operand.) 
For example, to do the sum of $"1.1 + 1.3"$, a computer converts these numbers into binary form and then does the math. Hence, because of the limitation in the number of digits it can store, the result is not exactly $2.4$; It is $2.4000000000000004$. 
\par It is known that floating-point addition and multiplication are both commutative but not necessarily associative nor distributive \cite{compsci1}. This type of error is the reason for the difference caused by using different mathematically equivalent expressions. In fact, this error is caused by a round-off error and is from the order of precision of the largest number in the operations.
\end{enumerate}
\par
Both of those errors mentioned above are inevitable, but it is important to recognize which one is larger and therefore is the main reason for the final error. \par
In our simulations, we used NumPy float32, which can store up to 7 digits, and the RK4 algorithm with discretization parameter $dt = 10^{-3}$. Since the numbers we encountered in calculations were at most $\mathcal{O}({10}^{2})$, the rounding errors are $\mathcal{O}({10}^{-5})$. We should compare this to the \emph{local truncation} error, which is $(10^{-3} )^5= 10^{-15}$; So, we can conclude that in this computation, RK4's error is negligible compared to that of round-off.
\section{Validity of the simulation result}
\label{sec:4}
\subsection{Method to determine validity of the numerical solution \& obtain a reliable result}\label{sec:4.1}
To discuss the validity of the numerical results, the first step is to find the main error of the simulation according to section \ref{sec:3}. The second step is to utilize the Lyapunov exponent to see how they affect the result. 
\par In transient chaos, before reaching the final non-chaotic destiny, the system undergoes a chaotic region and thus has a positive Lyapunov exponent \cite{lyap}. If the time the system behaves chaotically and its Lyapunov exponent is respectively denoted by $\Delta t$ and $\lambda$, then the deviation from the true solution will be:
\begin{equation}
    \delta \sim \delta_{0} e^{\lambda \Delta t} \label{eq:2}
\end{equation}
, where $\delta_0$ is the numerical error obtained from section \ref{sec:3}.
\par
If $\delta$ (the final error of the computation process) is more than the order of the system’s sensitivity, $S$, which is defined in \eqref{eq:definition}, it can be concluded that the result is wrong. However, it is not an easy task to calculate the sensitivity of the system. Therefore, we compare $\delta$ with the size of the strange attractor, which is an upper limit for $S$. If it is more than that, then the result is definitely not reliable. It is worth mentioning that $\delta$ needs to be much less than the size of the strange attractor for a reliable result. 

\begin{align}
        for\ \forall \epsilon<S:     |\overrightarrow{r}_{numerical}(0)-\overrightarrow{r}_{real}(0)|=\epsilon \notag
        \\
         \implies \lim\limits_{t \to \infty}    |\overrightarrow{r}_{numerical}(t)-\overrightarrow{r}_{real}(t)|<\epsilon     \label{eq:definition}
\end{align}

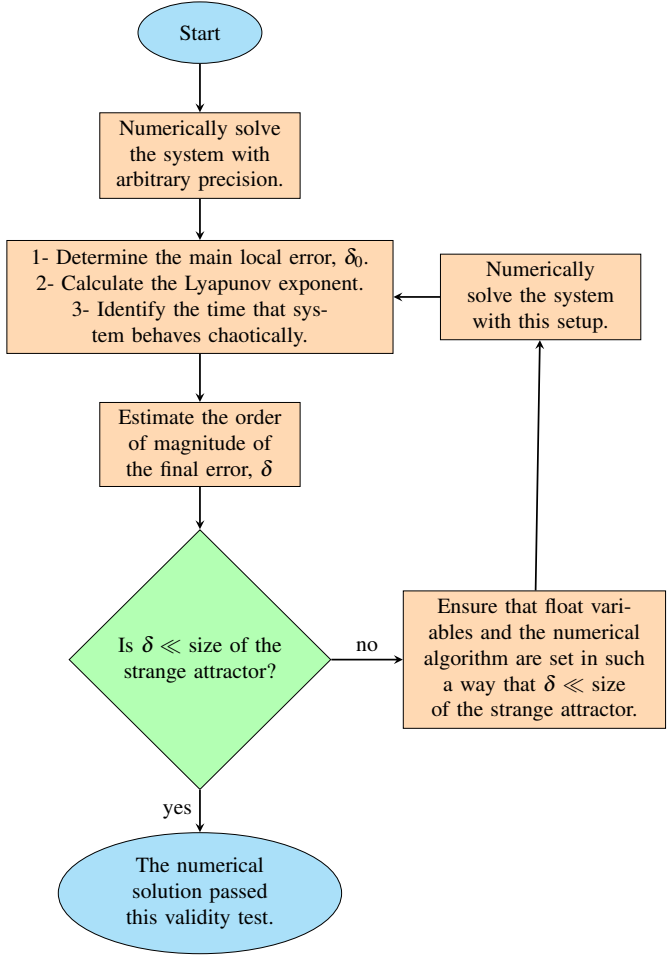
\begin{figure}[ht]
\centering 
\resizebox{0.5\textwidth}{!}{
\begin{tikzpicture}[node distance=2cm]

\node (start) [start] {Start};
\node (pro1) [process, below of=start] {Numerically solve the system with arbitrary precision.};
\node (pro2) [lprocess, below of=pro1, yshift=-0.3cm] {1- Determine the main local error, $\delta_0$. \\ 2- Calculate the Lyapunov exponent.\\ 3- Identify the time that system behaves chaotically.};
\node (pro3) [process, below of=pro2, yshift=-0.4cm] {Estimate the order of magnitude of the final error, $\delta$};

\node (dec1) [decision, below of=pro3, yshift=-1.5cm] {Is $\delta \ll$ size of the strange attractor?};

\node (pro4) [process, right of= pro2, xshift=3.5cm] {Numerically solve the system with this setup.};

\node (pro3b) [mprocess, right of=dec1, xshift=3.4cm] {Ensure that float variables and the numerical algorithm are set in such a way that $\delta \ll$ size of the strange attractor.};
\node (stop) [stop, below of=dec1, yshift= -1.8 cm] {The numerical solution passed this validity test.};

\draw [arrow] (start) -- (pro1);
\draw [arrow] (pro1) -- (pro2);
\draw [arrow] (pro2) -- (pro3);
\draw [arrow] (pro3) -- (dec1);
\draw [arrow] (dec1) -- node[anchor=east] {yes} (stop);
\draw [arrow] (dec1) -- node[anchor=south] {no} (pro3b);
\draw [arrow] (pro3b) -- (pro4);
\draw [arrow] (pro4) -- (pro2);

\end{tikzpicture}}

\caption{
The algorithm to determine reliability of the numerical solution and reach the actual result of the system.}
\label{fig:Flowchart}
\end{figure}


\begin{figure*}[!htp]
    \begin{subfigure}{0.45\linewidth}
    \includegraphics[width=\linewidth]{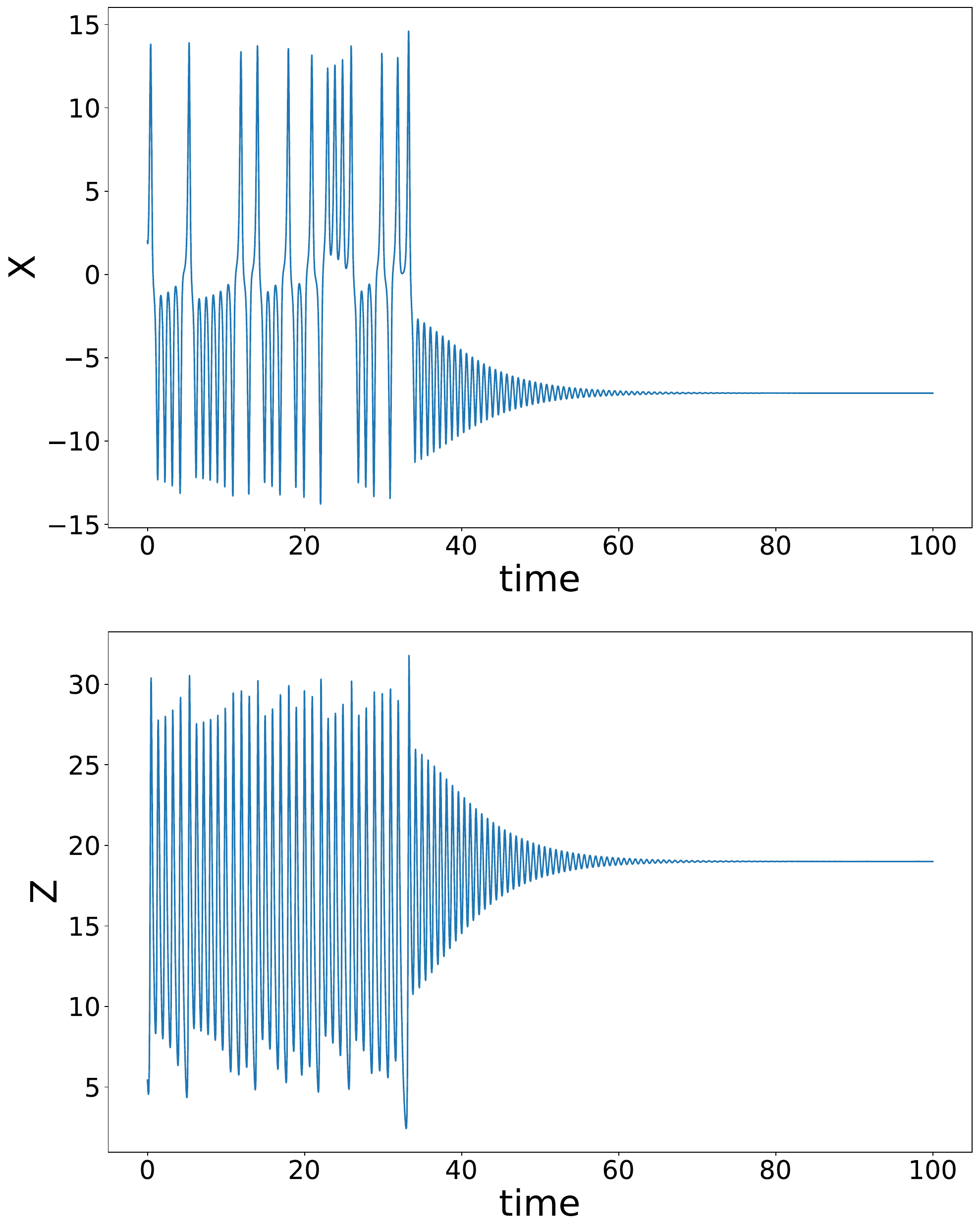} 
    \caption{with expression \eqref{eq:1.a}.}
    \label{fig:4.a}
    \end{subfigure}
    \hspace{0.2 mm}
    \begin{subfigure}{0.45\linewidth}
    \includegraphics[width=\linewidth]{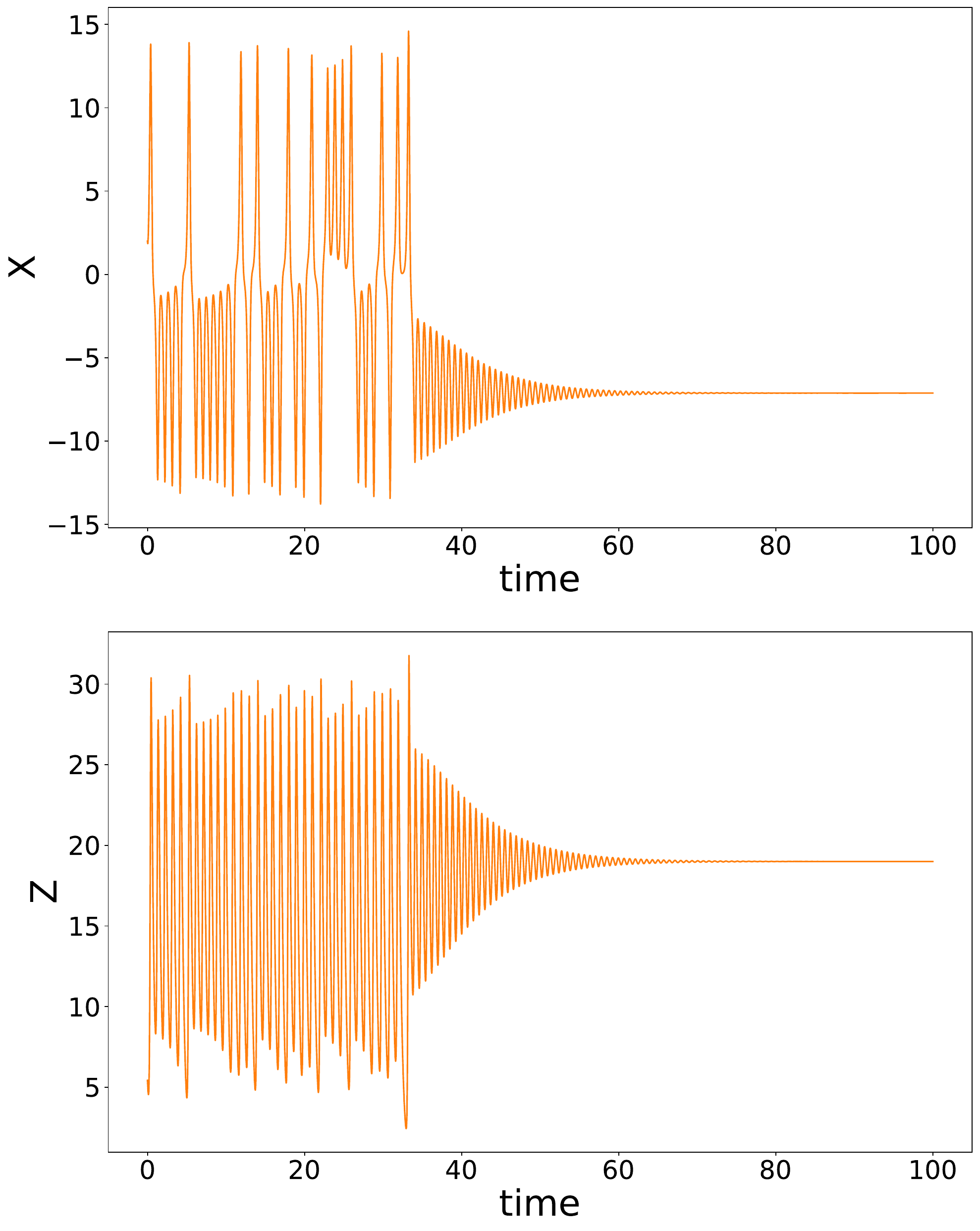}
    \caption{with expression \eqref{eq:1.b}.}
    \label{fig:4.b}
    \end{subfigure}
    \caption{Mathematically equivalent expressions reaching the same fixed point. \\Numerical results of the provided example with \textbf{NumPy float64} variables, X(t) \& Z(t) plots, where colors blue \& orange, \ref{fig:4.a} \& \ref{fig:4.b}, refer to the equations \eqref{eq:1.a} \& \eqref{eq:1.b} respectively: Lorenz equations with $r=20$ \& initial condition $(x_0,y_0,z_0)=(2,1,5.42857)$, solved with RK4 algorithm. \\The expressions are mathematically equivalent, and the system reaches the same destination, i.e. the fixed point $C_-$.}
    \label{fig:4}
\end{figure*}
\par
Suppose it was concluded that the result was wrong. In that case, one should reduce the order of main error by either changing the numerical method or setting float variables to increase the precision according to the type of main error (as mentioned in section \ref{sec:3}.) Thus, by reaching a precision where $\delta \ll$ size of the strange attractor, it can be concluded that the result might be valid. The flowchart of the proposed method is shown in Fig.\ref{fig:Flowchart}.

\subsection{Applying the method to the studied case}
\label{sec:4.2}

To apply the method proposed in section \ref{sec:4.1} to the provided example in section \ref{sec:2}, first we have to determine the main error. According to section \ref{sec:3}, $\delta_0 \sim 10^{-5}$ (round-off error) and from Fig.\ref{fig:2}, $\Delta t \approx 38$. The Lyapunov exponent of this system is $\lambda \approx 0.83$. Thus: 
\begin{equation*}
    \mathcal{O}(\delta)= \mathcal{O}(10^{-5}\times e^{0.83\times 38})=10^9
\end{equation*}
However, the size of the strange attractor in our system is $\mathcal{O}(10^2)$, and we know that $S$ cannot exceed it. As $10^2 \ll 10^9$, this indicates that our numerical solutions are completely wrong and explains why we encountered two different results using mathematically equivalent expressions as in Fig.\ref{fig:2}.
Thus, to get a reliable result, as stated in \ref{sec:4.1}, we need to increase the precision of the variables. So, in our simulation, we used NumPy float64 variables instead of NumPy float32, which can store up to 15-16 digits. Thus, as mentioned in section \ref{sec:3}, the round-off error in our system is $\mathcal{O}({10}^{-14})$. Now, we can estimate $\delta$ from \eqref{eq:2}:
\begin{equation*}
    \mathcal{O}(\delta)= \mathcal{O}(10^{-14}\times e^{0.83\times 38})=1
\end{equation*}
which is much less than the strange attractor’s size. But we cannot still be sure that float64 will result in a valid outcome; because we are not able to know $S$ precisely. Fig.\ref{fig:4} is produced by repeating our simulation exactly as in section \ref{sec:2}, with the only difference that we use NumPy float64 instead of NumPy float32. As it is seen, the results are so close, and the system reaches the same destination, i.e., fixed point, for expressions \eqref{eq:1.a} and \eqref{eq:1.b}. Thus, it is more probable that our simulation with NumPy float64 is correct. 
\par
However, it is of great importance to mention that getting the same results from mathematically equivalent expressions is only a necessary condition and not sufficient for the validation of our numerical solution.
\par
To ensure the reason behind obtaining different results for mathematically equivalent expressions is the rounding error, we reran the simulations of this article in C language, achieving the same results. This was expected since C compilers (specifically GCC) use the same floating-point arithmetic standard (IEEE 754) as Python and NumPy, confirming our claim.

\section{Summary}
In examining systems that undergo transient chaos, when the final destiny of the system is important, one should question reliability of the solution. The first step is to estimate the main error, as shown in section \ref{sec:3}. Then, by calculating the Lyapunov exponent and the time the system behaves chaotically, the final numerical error will be obtained from the equation \eqref{eq:2}. If the final error was more than the sensitivity of the system, as stated in section \ref{sec:4}, it could be concluded that the simulation result is wrong. If so, to get a reliable result, one can decrease the main error and redo the stated procedure.
\par
Finally, it should be mentioned that although we used the Lorenz system to demonstrate the proposed method, the arguments are general, and it is expected to work for any system exhibiting transient chaos.

\insertbibliography{References}

\end{document}